\documentclass[12pt,A4]{article}

\usepackage{graphics} 
\usepackage{epsfig} 
\usepackage{amsmath} 
\usepackage{amssymb}  
\usepackage{amscd}
\usepackage[latin1]{inputenc}
\usepackage{graphicx}

\newtheorem{thm}{Theorem}
\newtheorem{pf}{Proof}
\newtheorem{defn}{Definition}

\newtheorem{lem}{Lemma}
\newtheorem{assum}{Assumption}

\title{Reliable Control for Parameterized
Interconnected Systems using LMI techniques }

\author{Gisela Pujol, Jos\'{e} Rodellar, Josep M. Rossell
\thanks{This work was funded by CICYT-Spain under project DPI2002-04018-C02-01.
The partial support of the Government of Catalonia's grant
ACI053-2002 is also appreciated.}\\
Department of Applied Mathematics III\\
Technical University of Catalonia\\
Comte Urgell 187, 08036 Barcelona, Spain\\
\texttt{gisela.pujol@upc.edu, jose.rodellar@upc.edu,}\\
\texttt{josep.rossell@upc.edu} }

\begin{document}

\maketitle

\begin{abstract}
This work presents the design of a reliable decentralized state
feedback control  for a class of uncertain interconnected
polytopic continous systems. A model of failures in actuators is
adopted which considers outages or partial degradation in
independent actuators. The control is developed using the concept
of guaranteed cost control and a new LMI characterization using
polytopic Lyapunov functions. \textit{subj-class: OC -
Optimization and Control}. \textit{MSC-class:93A14,93B18,49K99}
\end{abstract}

\section{Introduction}

The decentralized control of interconnected systems in the
presence of uncertainties has attracted considerable attention in
the last years. With the aim of stabilizing the overall system
while ensuring a satisfactory performance, the guaranteed cost
control approach (GCC) has been recently considered (see
\cite{c3}, \cite{c4}, \cite{c6} and \cite{c7}). The guaranteed
cost control is concerned with the design of a state feedback
controller such that the closed-loop system is stable and an upper
bound of a quadratic cost function is minimized.

Another important issue when dealing with interconnected systems
is the design of fault-tolerant control systems. Reliable control
is concerned with the design of a closed-loop system to maintain
key properties, in spite of sensor or actuator outage or partial
degradation. Two main approaches have been proposed in the
literature. One uses multiple controllers in a redundant control
scheme (\cite{Sil80}, \cite{Yang98}). The other approach seeks a
reliable control design without redundancy by ensuring stability
and some performance bounds for specified class of admissible
failures of particular control components \cite{Vei95}. Within
this approach, the class of failures have been usually modelled as
outages.  This model considers that the control set can be
partitioned into two subsets: one subset includes the actuators
whose failures are admissible in the control design; and the
complementary subset with the actuators that are assumed to keep a
normal operation (see \cite{c9}). A more general failure scheme is
considered in \cite{c7}, where the authors present a model of
actuator failures which takes into account the outage case and
also the possibility of partial failures. In \cite{c7} a
centralized reliable control is designed for a class of  uncertain
nonlinear systems without subsystem interconnections. In the
present paper, we adopt this failure model extended within a
reliable decentralized control scheme for uncertain interconnected
linear systems. Linear matrix inequalities (LMI) techniques have
been used in the context of guaranteed cost control problems and
reliable control design (\cite{Yan01}, \cite{c9}). The currently 
known LMI characterizations are potentially conservative in the
sense that they use a common Lyapunov function regardless of the
parameter values. In order to reduce conservatism in the case of
constant parameters, the notion of parameter-dependent function
was introduced \cite{c2}. The Lyapunov conditions leads to
nonconvex optimization problem which does not seem tractable in
general. In \cite{c1} and \cite{c5}, this weakness is overcome
using auxiliary variables, as well as replacing the functions of
singular Lyapunov for multiple functions  to obtain more robust
tools and to reduce the conservatism in the control problems, and
provides additional flexibility in a wide range of problems. In
\cite{c5}, a way to "convexifying" the general affine Lyapunov
problem has been proposed, leading to a parametric linear matrix
inequality.

This paper has two objectives. One is related to the decentralized
reliable guaranteed cost control problem for interconnected
continous systems. To our knowledge, there are not available
results dealing with reliable control for uncertain interconnected
systems with guaranteed cost. We present a control design that
allows to consider reliability and guaranteed cost. This reliable
control shows that the admission of control failures imposes some
restriction in the control weighting matrices in the performance
criterion. Thus the designer can take some trade-off between
control performance and admitted reliability. The second objective
is to give a new LMI characterization and its relation with
polytopic systems and polytopic Lyapunov functions. In the
centralized control design of uncertain systems with guaranteed
cost and reability, matrix inequalities are obtained which lead to
LMI's through a linealization process. When dealing with
decentralized control of interconnected systems this process needs
to be extended. In this paper, we obtain a general result that can
be applied to different classes of problems, such as guaranteed
cost, $H_{2}$ or $H_{\infty}$ control design. We apply this result
in the design of decentralized reliable control of uncertain
interconnected systems with guaranteed performance.

The notation throughout the paper is fairly standard. In symmetric
block matrices or long matrix expressions, we use $*$ as an
ellipsis for terms that are induced by symmetry, e.g.
$$
\left( \begin{array}{cc} S+(*) & *\\ M&Q \end{array} \right)\equiv
\left( \begin{array}{cc} S+S^{T} & M^{T}\\ M&Q
\end{array} \right) \, .
$$
We use bold to denote dependence on the parameter $\alpha$, e.g.
$\textbf{A}=A(\alpha)$; $I$ denotes the identity matrix.

\section{Problem statement}
Consider a class of large-scale interconnected system composed of
$N$ polytopic subsystems described by the following state
equations:
\begin{equation}
\label{sis1} \left\{
\begin{array}{l} \dot{x}_{i}=A_{i}(\alpha)x_{i} + B_{i}(\alpha)u_{i}+\sum_{j
\neq i} G_{ij} g_{ij}(t,x_{j}) \\[1.5ex]
x_{i}(0)=x_{i0} , \; i=1,\dots,N\;.
\end{array} \right.
\end{equation}
\noindent We define:
\begin{itemize}
\item $x_{i} \in \mathbb{R}^{n_{i}}$ state of the $i$th subsystem
\item $u_{i} \in \mathbb{R}^{s_{i}}$ control of the $i$th
subsystem \item $A_{i}(\alpha) \in \mathbb{R}^{n_{i} \times
n_{i}}$ parametric state matrix \item $B_{i}(\alpha) \in
\mathbb{R} ^{n_{i} \times s_{i}}$ parametric control matrix \item
$g_{ij}(t,x_{j}) \in \mathbb{R}^{l_{i}}$ unknown interconnection
vector function
 \item $G_{ij} \in \mathbb{R}^{n_{i}\times l_{i}}$ constant
interconnection matrix
\end{itemize}
It is assumed that the unknown vectors $g_{ij}(t,x_{j})$ are
continuous and sufficiently smooth in $x_{j}$ and piecewise
continuous in $t$.

\textit{Uncertain model }

\vspace{-0.2cm} The parameter uncertainties considered here are
assumed to be of the following form:
\begin{equation}
\label{poli}
A_{i}(\alpha)=\sum_{k=1}^{L} \alpha_{k}A_{i_{k}} \;,\;
B_{i}(\alpha)=\sum_{k=1}^{L} \alpha_{k}B_{i_{k}}\;.
\end{equation}
\noindent The parameter vector $\alpha$ is in the simplex $\Pi$
defined by
\begin{equation}
\label{polipi} \Pi=\{\alpha \in \mathbb{R}^{L},\;\sum_{k=1}^{L}
\alpha_{k} =1,\; \alpha_{k} \geq 0,\; k=1,\dots, L \}\;.
\end{equation}
\noindent That is, $A_{i}(\alpha)\in
\textrm{co}\{A_{i_{1}},\dots,A_{i_{L}} \}$ and $B_{i}(\alpha)\in
\textrm{co}\{B_{i_{1}},$ $\dots,B_{i_{L}} \}$ are convex
combinations of the matrix $A_{k}$ and $B_{k}$, respectively.

\textit{Interconnection assumptions} \vspace{-0.2cm}
\begin{assum}
\label{asum1} There exist known constant matrices $ W_{ij}$ such
that, for all $x_{j} \in R^{n_{j}}$,
\begin{equation}
\label{gW} || g_{ij}(t,x_{j}) || \leq || W_{ij} \; x_{j} ||
\end{equation}
for all $i$, $j$ and for all $t \geq 0$, where $||\;||$ denotes
the Euclidean norm.
\end{assum}
\begin{assum}
\label{asum2} For all $i$, $\small{W_{i}:=\sum\limits_{j=1,j \neq
i}^{N}W_{ji}^{T} \; W_{ji}>0}\;.$
\end{assum}\vspace{-0.4cm}
The Assumption \ref{asum1} allows some linear structure to the
 interconnection vector $g_{ij}$ in terms of the state vector
$x_{j}$. The Assumption \ref{asum2} ensure that almost two
subsystems must be interconnected, because at least one $W_{ij}$
must be nonzero, for each $i$.

\textit{Cost function}

\vspace{-0.2cm}Consider the following cost function associated
with system (\ref{sis1}):
\begin{equation}
\label{cost} J(x,u)=\sum_{i=1}^{N} \int_{0}^{\infty}
\left(x_{i}^{T}Q_{i}x_{i} + u_{i}^{T}R_{i}u_{i}\right)\, dt
\end{equation}
\noindent where $Q_{i} \in \mathbb{R}^{n_{i} \times n_{i}}$ and
$R_{i} \in \mathbb{R}^{s_{i} \times s_{i}}$ are given constant
symmetric positive definite  matrices, for all $i$.

\textit{Failure model}

\vspace{-0.2cm}Let $u_{i}^{F}$ denote the vectors with the signals
from the $s_i$ actuators which control the $ith$ subsystem. Here
we consider the following failure model:
\begin{equation}
\label{fial} u_{i}^{F}=\Lambda_{i}u_{i}+\phi_{i}(u_{i}), \qquad
i=1,\dots,N
\end{equation}
\noindent where $\Lambda_{i}$=diag($\lambda_{i1},\dots,\lambda_{i
s_{i}}$)$\in \mathbb{R}^{s_{i} \times s_{i}}$ is a diagonal
positive definite matrix. The uncertain function
$\phi_{i}(u_{i})=(\phi_{i1}(u_{i1}),\dots,\phi_{is_i}(u_{i\,s_{i}}))$
satisfies, for each $i$,
\begin{equation}
\label{compfi} \phi_{ij}^{2}(u_{ij}) \leq
\gamma_{ij}^{2}\;u_{ij}^{2}\; , \; j=1,\dots,s_{i}
\end{equation}
\noindent where $\gamma_{ij} \geq 0$. If (\ref{compfi}) holds,
then
\begin{equation}
\label{fi} ||\phi_{i}(u_{i})||^{2} \leq ||\Gamma_{i}u_{i} ||^{2},
\qquad i=1,\dots,N
\end{equation}
\noindent where $\Gamma_{i} =$diag($\gamma_{i1},\dots,\gamma_{i
s_{i}}$)$\in \mathbb{R}^{s_{i} \times s_{i}}$ is a diagonal
positive semidefinite matrix.

The value of $\lambda_{ij}$ represents the percentage of failure
in the $j$th actuator of the controller of the $i$th subsystem.
Each actuator can fail independently. If $\lambda_{ij}=1$ and
$\gamma_{ij}=0$, it corresponds to the normal case for the $j$th
actuator of the $i$th subsystem ($u_{ij}^{F}=u_{ij}$). When this
is true for all $j$, we have $\Lambda_{i}=I_{s_{i}}$ and
$\Gamma_{i}=0$ and it corresponds to the normal case in the $i$th
canal ($u_{i}^{F}=u_{i}$). When $\lambda_{ij}=\gamma_{ij}$,
(\ref{fial}) and (\ref{compfi}) cover the outage case
($u_{ij}^{F}=0$) because $\phi_{ij}=-\lambda_{ij}\;u_{ij}$ verify
(\ref{compfi}). The case $\phi_{i}(u_{i})=-\Lambda_{i}u_{i}$
corresponds to the outage of the whole controller of the $i$th
system. Other cases correspond to partial failures or partial
degradations of the actuators.

\textit{Control objective (Reliable guaranteed cost
control)}

\vspace{-0.2cm}The objective of this paper is to design a set of
decentralized feedback control laws $u_{i}(t) = K_{i}x_{i}(t)$
 ($i=1,\dots,N$) and obtain a Lyapunov function defined by $\textbf{X}_i$,
  for the interconnected systems (\ref{sis1}) with
uncertainties model (\ref{poli})-(\ref{polipi}) and Assumptions 1
and 2, in such a way that, in the presence of the failures
described by (\ref{fial}) and (\ref{compfi}), the following
property is satisfied:
\begin{equation}
\label{qgcc} \sum_{i=1}^{N}\Big(
\frac{d}{dt}x_{i}^{T}\textbf{X}_{i}x_{i}+ x_{i}^{T}
Q_{i}x_{i}+(u_{i}^{F})^{T}R_{i}u_{i}^{F} \Big)<0\;.
\end{equation}
\noindent This inequality leads to a bound for the cost function
(\ref{cost}) in the form $J(x,u^{F}) \leq \bar{J}$, where
$\bar{J}$ is some specified constant, and to ensure that the
closed-loop system
\begin{equation}
\label{sis3} \left\{ \begin{array}{l} \dot{x}_{i}=A_{i}(\alpha)
x_{i}+ B_{i}(\alpha)u_{i}^{F} + \sum_{j=1,j \neq i}^{N}
G_{ij} \; g_{ij}\\[2ex]
u_{i}^{F}=\Lambda_{i}u_{i}+\phi_{i}(u_{i})\\[2ex]
u_{i}=K_{i}x_{i}
\end{array} \right.
\end{equation}
\noindent is asymptotically stable.
\begin{defn}
\label{rgcc} The set of above feedback control laws $u_{i}$ is
said to be a {\it reliable guaranteed cost control}.
\end{defn}


\section{Main results}

\subsection{Instrumental tools}
\begin{lem}[Projection lemma]
\label{lem1} For a symmetric matrix $\Psi$ and matrices $P$, $Q$
with appropriate dimensions, there exists a matrix $X$ such that:
$$ \Psi + P^{T}X^{T}Q+Q^{T}XP<0 \Leftrightarrow
\left\{\begin{array}{l} \mathcal{N}_{P}^{T}\Psi
\mathcal{N}_{P}<0 \\[2ex]
\mathcal{N}_{Q}^{T}\Psi \mathcal{N}_{Q}<0
\end{array}\right.
$$
\noindent with $\mathcal{N}_{P}$ and $\mathcal{N}_{Q}$ any
matrices whose columns form bases of  $P$ and $Q$ respectively.
\end{lem}
\begin{lem}[Schur Complement]
\label{schur} Consider a symmetric matrix
$M=\left(\begin{array}{cc}P_{1}&P_{2}\\P_{2}^{T}&P_{3}
\end{array} \right)$. Then,
$$
M<0 \Leftrightarrow \left\{\begin{array}{l}
P_{1}-P_{2}P_{3}^{-1}P_{2}^{T}<0 \\ [1ex] P_{3}<0 \qquad
invertible.
\end{array} \right.
$$
\end{lem}
\subsection{LMI characterization}
Here we give a general theorem that allows a LMI characterization.
The subindex $i$ is omitted because the results can be applied to
interconnected systems but also to a single system.
\begin{thm}
\label{tem1} The following statements, involving symmetric
positive definite  matrix variables $\textbf{X}$, $\textbf{Y}$ and
general matrix variables $V$, $N$, $K$, are equivalent.
\begin{enumerate}
\item[(i)] There exist $\textbf{X}$, $K$ such that
{\footnotesize\begin{equation} \label{t11} \left(
\begin{array}{ccc} (\textbf{A} +
\textbf{B}K )^{T}\textbf{X} + (*)&*&*\\
G^{T}\textbf{X}&Q_{11}&*\\
C&Q_{21}&Q_{22}
\end{array} \right)<0
\end{equation}}

\item[(ii)] There exist $\textbf{Y}$, $V$ and $N$ such that
{\footnotesize\begin{equation} \label{t12} \left(
\begin{array}{ccccc}
-(V+V^{T})&*&*&*&*\\
\textbf{A}V+\textbf{Y}+\textbf{B}N& -\textbf{Y}&*&*&*\\
0&G^{T}&Q_{11}&*&*\\
CV&0&Q_{21}&Q_{22}&*\\
V&0&0&0&-\textbf{Y}
\end{array}
\right)<0
\end{equation}
\begin{equation}
\label{t122} \left( \begin{array}{ccc}
-\textbf{Y}&*&*\\
G^{T}&Q_{11}&*\\
0&Q_{21}&Q_{22}
\end{array}
\right)<0
\end{equation}}

\item[(iii)] Consider $\textbf{A}=A(\alpha)$,
$\textbf{B}=B(\alpha)$, $\textbf{X}=X(\alpha)$ as in (\ref{poli}).
There exist $\textbf{Y}=Y(\alpha)$, $V$ and $N$ such that
{\footnotesize\begin{equation} \label{t13} \left(
\begin{array}{ccccc}
-(V+V^{T})&*&*&*&*\\
A_k V+Y_k+B_k N&-Y_k&*&*&*\\
0&G^{T}&Q_{11}&*&*\\
CV&0&Q_{21}&Q_{22}&*\\
V&0&0&0&-Y_k
\end{array}
\right)<0
\end{equation}
\begin{equation}
\label{t132} \left( \begin{array}{ccc}
-Y_k&*&*\\
G^{T}&Q_{11}&*\\
0&Q_{21}&Q_{22}
\end{array}
\right)<0\;.
\end{equation}}
\end{enumerate}
\end{thm}
\vspace*{0.2cm} The matrix $Q=\left(
\begin{array}{cc}Q_{11} &*\\Q_{21}&Q_{22}
\end{array} \right)$ has to be a negative defined matrix.
\vspace{-0.4cm}
\begin{pf}
Taking into account that $\sum_{k=1}^{L} \alpha_{k}=1$, the
inequality (\ref{t12}) is true when it holds in the vertices of
the simplex $\Pi$, so (ii) and (iii) are equivalent. To see the
equivalence between (i) and (ii), consider a new variable
$N_{i}:=K_{i}V_{i}$ and apply Projection Lemma \ref{lem1} to
(\ref{t12}). Consider $\bar{\textbf{A}}=\textbf{A}+\textbf{B}K$,
$P=\left( Id,0,0,0,0 \right)$, $Q=\left(
-Id,\;\bar{\textbf{A}}^{T},\;0,\;C^{T},\;Id \right)$ and
{\footnotesize $$\psi=\left( \begin{array}{ccccc}
0&*&*&*&*\\
\textbf{Y}&-\textbf{Y}&*&*&*\\
0&G^{T}&Q_{11}&*&*\\
0&0&Q_{21}&Q_{22}&*\\
0&0&0&0&-\textbf{Y}
\end{array}
\right)\; .$$} \noindent The null spaces bases of $P$ and $Q$ are
$${\footnotesize
\mathcal{N}_{P}=\left( \begin{array}{cccc}
0&0&0&0\\
I&0&0&0\\
0&I&0&0\\
0&0&I&0\\
0&0&0&I
\end{array}
\right), \,\,\, \mathcal{N}_{Q}=\left( \begin{array}{cccc}
\bar{\textbf{A}}^{T}&0&C^{T}&Id\\
I&0&0&0\\
0&I&0&0\\
0&0&I&0\\
0&0&0&I
\end{array}
\right)}\;.
$$
\noindent By $\mathcal{N}_{P}^{T}\psi\mathcal{N}_{P}<0$, we have
{\footnotesize\begin{equation} \label{liap31} \left(
\begin{array}{cccc}
-\textbf{Y}&*&*&*\\
G^{T}&Q_{11}&*&*\\
0&Q_{21}&Q_{22}&*\\
0&0&0&-\textbf{Y}
\end{array}
\right)<0.
\end{equation}}
\noindent By Schur complement, (\ref{liap31}) is equivalent to
(\ref{t122}).

From $\mathcal{N}_{Q}^{T}\psi\mathcal{N}_{Q}<0$, we have
{\footnotesize$$ \left( \begin{array}{cccc}
\bar{\textbf{A}}\textbf{Y}+\textbf{Y}\bar{\textbf{A}}^{T}- \textbf{Y}&*&*&*\\
G^{T}&Q_{11}&*&*\\
C\textbf{Y}&Q_{21}&Q_{22}&*\\
\textbf{Y}&0&0&-\textbf{Y}
\end{array}
\right)<0\;.
$$}
\noindent By Schur complement, {\footnotesize\begin{equation}
\label{11} \left( \begin{array}{ccc}
\bar{\textbf{A}}\textbf{Y}+\textbf{Y}\bar{\textbf{A}}^{T}&*&*\\
G^T&Q_{11}&*\\
C\textbf{Y}&Q_{21}&Q_{22}
\end{array}
\right)<0\;.
\end{equation}}
\noindent Consider $Q=\left(\begin{array}{cc}Q_{11}
&*\\Q_{21}&Q_{22}
\end{array} \right)$. Applying Schur complement again,
{\footnotesize$$
\bar{\textbf{A}}\textbf{Y}+\textbf{Y}\bar{\textbf{A}}^{T}-(G \quad
\textbf{Y}C^{T})Q^{-1} \left(\begin{array}{c}G^{T}\\C\textbf{Y}
\end{array} \right)<0\;.
$$}
\noindent Using now that $\textbf{Y}=\textbf{X}^{-1}$, we have
{\footnotesize$$
\textbf{X}\bar{\textbf{A}}+\bar{\textbf{A}}^{T}\textbf{X}-(\textbf{X}G
\quad C^{T})Q^{-1} \left(\begin{array}{c}G^{T}\textbf{X}\\C
\end{array} \right)<0,
$$}
\noindent which is equivalent to (\ref{t11}).
\hspace{3cm}$\square$\end{pf}\vspace{-0.3cm}

\subsection{Stability problem}
In this subsection, we apply Theorem \ref{tem1} to stability
characterization. We introduce an alternative characterization of
the fundamental Lyapunov stability theorem for linear
interconnected systems.  It introduces a new transformation on the
Lyapunov variables which helps to reduce the typical degree of
conservatism in some problems. The obtained result will be the
base of the development of the reliable guaranteed control in the
next section.

Consider the polytopic Lyapunov function:
\begin{equation}
\label{liap} V(x,\alpha) =
\sum_{i=1}^{N}x_{i}^{T}\textbf{X}_{i}x_{i},
\end{equation}
\noindent where $\textbf{X}_{i}=\sum_{k=1}^{L}
\alpha_{k}X_{i_{k}}$, $\alpha \in \Pi$ and $X_{i_{k}}>0$, for all
$k=1,\dots,L$ and  for all $i=1,\dots,N$. This function  has to
satisfy:
\begin{enumerate}
\item[(i)] $V(x,\alpha) > 0$, \item[(ii)] $\dot{V}(x,\alpha)<0$.
\end{enumerate}
Let us evaluate condition (ii). From (\ref{liap}),
$$
\dot{V}(x,\alpha)=\sum_{i=1}^{N}\dot{x}_{i}^{T}\textbf{X}_{i}x_{i}+x_{i}^{T}\textbf{X}_{i}\dot{x}_{i}\;<0.
$$
\noindent By using the feedback controls $u_{i}=K_{i}x_{i}$ in
(\ref{sis1}),
$$
\dot{V}=\sum_{i=1}^{N}x_{i}^{T}\big( (\textbf{A}_{i} +
\textbf{B}_{i})K_{i})^{T}\textbf{X}_{i} +
\textbf{X}_{i}(\textbf{A}_{i} + \textbf{B}_{i}K_{i} ) \big)x_{i}+
$$
$$
+ \sum_{i=1}^{N}\sum_{j \neq
i}g_{ij}^{T}(t,x_{j})\;G_{ij}^{T}\textbf{X}_{i}x_{i}
+x_{i}^{T}\textbf{X}_{i}G_{ij}\;g_{ij}(t,x_{j})<0\, .
$$
 \noindent Now, using (\ref{gW}), this inequality leads to
{\footnotesize\begin{equation} \label{8} \left( \begin{array}{ccc}
(\textbf{A}_{i} +
\textbf{B}_{i}K_{i} )^{T}\textbf{X}_{i} + (*)&*&*\\
G_{i}^{T}\textbf{X}_{i}&Q_{11}^{i}&*\\
C_{i}&Q_{21}^{i}&Q_{22}^{i}
\end{array} \right)<0
\end{equation}}
\noindent with $ G_{i}=(G_{i1},\dots,G_{iN})$, $C_{i} = Id$, $
Q_{11}^{i}=diag(-Id,$ $,\dots,-Id)$, $Q_{21}^{i}=0$, $
Q_{22}^{i}=-W_{i}^{-1}$. This step involves tedious manipulations.
The details are omitted here, but can be found in \cite{c8}. So
that, if (\ref{8}) holds for each $i=1,\dots,N$, then $\dot{V}<0$
and the system is stable.  In order to obtain an efficient LMI
condition, we now use Theorem \ref{tem1} to characterize the
stability.
\begin{thm}
\label{temest} Under assumptions \ref{asum1} and \ref{asum2},
consider the system (\ref{sis1}) with uncertain model (\ref{poli})
and (\ref{polipi}), and consider the state-feedback controls
$u_{i}=K_{i}x_{i}$. Assume that there exist symmetric positive
definite matrices $\left\{Y_{i_{k}}\right\}_{k}$ and matrices
$V_{i}$ and $N_{i}$ such that the LMI system
{\footnotesize\begin{equation} \label{10} \left(
\begin{array}{ccccc}
-(V_{i}+V_{i}^{T})&*&*&*&*\\
A_{i_{k}}V_{i}+Y_{i_{k}}+B_{i_{k}}N_{i}&- Y_{i_{k}}&*&*&*\\
0&G_{i}^{T}&Q_{11}^{i}&*&*\\
C_{i}V_{i}&0&Q_{21}^{i}&Q_{22}^{i}&*\\
V_{i}&0&0&0&-Y_{i_{k}}
\end{array}
\right)<0
\end{equation}}
{\footnotesize$$ \left( \begin{array}{cc}
-Y_{i_{k}}&*\\
G_{i}^{T}&Q_{11}^{i}
\end{array}
\right)<0
$$}
\noindent is feasible for all $k=1,\dots,L$ and $i=1,\dots,N$.
Then, the system is stable with Lyapunov function (\ref{liap}) and
{\footnotesize $$ \left\{
\begin{array}{l}
K_{i}=N_{i}V_{i}^{-1}\\[1.5ex]
X_{i}(\alpha)=\sum_{k=1}^{L}\alpha_{k}Y_{i_{k}}^{-1}
\end{array}
\right.
$$}
\end{thm}
\vspace{-0.4cm}\begin{pf}The proof is an immediate application of
Theorem \ref{tem1}, for each $k$.\hspace{3.27cm}
$\square$\end{pf}\vspace{-0.4cm} If we use a Lyapunov function
parameter independent, we can consider $Y_{i}$ instead of
$Y_{i_{k}}$ and the same expression is obtained.
\subsection{Reliable guaranteed cost control}
By Definition \ref{rgcc}, the system (\ref{sis3}) has reliable
guaranteed cost control $u_{i}=K_{i}x_{i}$ if there exist
symmetric variable matrices  $\textbf{X}_{i}$ and  gain matrices
$K_{i}$ such that the following inequality
{\footnotesize\begin{equation} \label{M2} \left(
\begin{array}{ccc} \Xi_{i}&*&* \\
G_{i}^{T}\textbf{X}_{i}&-I&*\\
\textbf{B}_{i}^{T}\textbf{X}_{i}+R_{i}\Lambda_{i}K_{i}&0&R_{i}-I
\end{array}
\right)<0,
\end{equation}}
\noindent is feasible, where $$\Xi_{i}=(\textbf{A}_{i} +
\textbf{B}_{i}\Lambda_{i}K_{i} )^{T}\textbf{X}_{i} +
\textbf{X}_{i}(\textbf{A}_{i} + \textbf{B}_{i}\Lambda_{i}K_{i}
)+$$$$+W_{i}+Q_{i}+K_{i}^{T}\Gamma_{i}^{2}K_{i}+K_{i}^{T}\Lambda_{i}R_{i}\Lambda_{i}K_{i}\;.$$\noindent
This characterization involves tedious manipulations where
(\ref{gW}) and (\ref{fi}) are used. The details are omitted here,
but can be found in \cite{c8}. By Schur complement, inequality
(\ref{M2}) is equivalent to (\ref{ineq}), defined in figure
\ref{ineq1}.
\begin{figure*}[hbt!]
\hrule
\vbox to 4.5cm{\vfill 
\begin{equation}\label{ineq}
{\small\left(
\begin{array}{ccccccc}
(\textbf{A}_{i}+\textbf{B}_{i}\Lambda_{i}K_{i})^{T}\textbf{X}_{i}+(*)&*&*&*&*&*&*\\[2ex]
G_{i}^{T}\textbf{X}_{i}&-I&*&*&*&*&* \\[2ex]
\textbf{B}_{i}^{T}\textbf{X}_{i}+R_{i}\Lambda_{i}K_{i}&0&R_{i}-I&*&*&*&*\\[2ex]
\Lambda_{i}K&0&0&-R_{i}^{-1}&*&*&*\\[2ex]
\Gamma_{i}K_{i}&0&0&0&-I&*&*\\[2ex]
I&0&0&0& 0&-Q_{i}^{-1}&*\\[2ex]
I&0&0&0&0&0&-W_{i}^{-1}\\[2ex]
\end{array}
\right)<0\;.}
\end{equation}
\vfill} \hrule\caption{}\label{ineq1}
\end{figure*}

In order to obtain an LMI characterization, it is necessary to
separate the terms in $\textbf{X}_{i}$ and $K_{i}$ in $\Xi_i$ from
(\ref{M2}) to be linear. After some manipulations, we obtain that
this is possible if $R_{i}-Id$ is a negative definite matrix and
invertible. In this way, we obtain a restriction on the cost
function. This means that we may loose some freedom in prescribing
the control performance to achieve a reliable control. In
accordance with this result, we introduce the following assumption
before constructing the reliable guaranteed cost control.

\textit{Assumption 3.} The cost control matrix $R_{i}$ must be
invertible and verify $R_{i}-Id<0$.\\

Taking this assumption into account, we obtain that (\ref{ineq})
is equivalent to: {\footnotesize\begin{equation} \label{M3} \left(
\begin{array}{cccc} (\textbf{A}_{i} + \hat{\textbf{B}}_{i}K_{i})^{T}\textbf{X}_{i}+(*)&*&*&* \\[2ex]
\textbf{E}_{i}^{T}\textbf{X}_{i}&Q_{11}^{i}&*&*\\[2ex]
C_{i}&0&Q_{22}^{i}&*\\[2ex]
F_{i}K_{i}&0&0&Q_{33}^{i}
\end{array}
\right)<0
\end{equation}}
\noindent where {\footnotesize$$ \hat{\textbf{B}}_{i}=
\textbf{B}_{i}(I+(Id-R_{i})^{-1}R_{i})\Lambda_{i}$$
$$
F_{i}=(\Lambda_{i},\Gamma_{i},R_{i}\Lambda_{i})\; , \;
\textbf{E}_{i}=(G_{i},\textbf{B}_{i})\;,\;C_{i}=(I,I)^{T}
$$
$$
Q_{11}^{i}=\left(\begin{array}{cc}-I&*\\0&R_{i}-I\end{array}\right),\;
\;
Q_{22}^{i}=\left(\begin{array}{cc}-Q_{i}^{-1}&*\\0&-W_{i}^{-1}\end{array}\right)$$
$$
Q_{33}^{i}=\left(\begin{array}{ccc}-R_{i}^{-1}&*&*\\0&-I&*\\0&0&R_{i}-I\end{array}\right).
$$}
\noindent We present now the main result.
\begin{thm}
\label{temestfiab} Under Assumptions \ref{asum1}, \ref{asum2} and
3, consider the system (\ref{sis3}) with the polytopic structure  
(\ref{poli}) and (\ref{polipi}). Suppose that, for each
$i=1,\dots,N$, there exist symmetric  positive definite matrices
$\left\{Y_{i_{k}}\right\}_{k=1,\dots,L}$ and matrices $V_{i}$ and
$N_{i}$  that  make the LMI system
\begin{equation}
\label{M1} {\footnotesize \left( \begin{array}{cccccc}
-(V_{i}+V_{i}^{T})&*&*&*&*&*\\
A_{i_{k}}V_{i}+Y_{i_{k}}+\hat{B}_{i_{k}} N_{i}&-Y_{i_{k}}&*&*&*&*\\
0&E_{i_{k}}^{T}&Q_{11}^{i}&*&*&*\\
C_{i}V_{i}&0&0&Q_{22}^{i}&*&*\\
F_{i}^{T}N_{i}&0&0&0&Q_{33}^{i}&*\\
V_{i}&0&0&0&0&-Y_{i_{k}}
\end{array}
\right)<0 }
\end{equation}
$$
\left(
\begin{array}{cc}
-Y_{i_{k}}&*\\
E_{i_{k}}^{T}&Q_{11}^{i}
\end{array}
\right)<0\;
$$
\noindent feasible for all $k=1,\dots,L$ and all $i=1,\dots,N$.
Then, the set of state feedback controls $u_{i}=K_{i}x_{i}$ is a
reliable guaranteed cost control with Lyapunov function
$V(x,\alpha)=\sum_{i=1}^{N}x_{i}^{T}X_{i}(\alpha)x_{i}$, where
\begin{equation}
\label{assig}  \left\{
\begin{array}{l}
K_{i}=N_{i}V_{i}^{-1}\\[2ex]
X_{i}(\alpha)=\sum_{k=1}^{L}\alpha_{k}Y_{i_{k}}^{-1}.
\end{array}
\right.\;
\end{equation}
\noindent Moreover, for any $x_{i}(0)=x_{i0}$, the cost function
satisfies
\begin{equation}\label{cotafc} J <
\sum_{i=1}^{N}x_{i0}^{T}\;X_{i}(\alpha)\;x_{i0}\; .
\end{equation}
\end{thm}
\vspace{-0.4cm}\begin{pf} Applying the Theorem \ref{tem1}, we
obtain (\ref{M3}) from (\ref{M1}).$\hspace{5cm}\square$
\end{pf}\vspace{-0.4cm}
The cost bound (\ref{cotafc}) depends on the initial conditions.
In order to eliminate this dependence, the mean value of the cost
function is sought over all possible values $x_{i0}$. This is
equivalent to:
$$
\mathcal{E}(J)<\sum_{i=1}^{N}tr(X_{i}(\alpha))=\sum_{i=1}^{N}\sum_{k=1}^{L}tr(X_{i_{k}})\alpha_{k},
$$
\noindent where $tr$ denotes the trace of a matrix. This equation
allows to find an optimum value. Considering
$\bar{J}:=\mathrm{max}_{k}\left(\sum_{i=1}^{N}tr(X_{i_{k}})\right)$,
we have the cost function bounded by $\mathcal{E}(J)<\bar{J}$.

\vspace{-3mm}

\section{Conclusions}

In this work, a solution for a reliable decentralized guaranteed
cost control problem for interconnected systems with  polytopic
uncertainties has been presented. Failures are described by a
model which considers possible outage or partial failures in every
actuator of each decentralized controller. The control design
involves two steps. First, an LMI characterization is presented.
Second, a sufficient condition is given for the existence of a
decentralized reliable guaranteed cost control set. A key point in
the control design has been the formulation of a new LMI
characterization, which uses parameter-dependent Lyapunov
functions and slack variables. The obtained LMI separates the
unknown variables from the system parameter data, which smoothes
the numerical solution. This characterization can be useful for
different class of problems, such as guaranteed cost control,
$H_2$ or $H_\infty$ control design. In the paper, this type of LMI
has been exploited to proof that the proposed decentralized
control scheme guarantees the quadratic stability and a cost bound
for a class of failure model which considers outage or partial
degradation of any independent specific actuator. The design
presented in this paper shows that the admission of control
failures imposes some restrictions in the definition of the cost
function to be bounded. Specifically, some freedom is lost in the
selection of the control weighting matrices.


\begin{thebibliography}{99}

\bibitem{c1}
P. Apkarian, H.D. Tuan and J. Bernussou, "Continous-Time Analysis,
eigenstructure Assignment and $H_{2}$ synthesis with Enhanced LMI
Characterizations", {\it in Proc. of 39th Conference on Decision
and Control}, Sydney, pp. 1489--1494, 2000.

\bibitem{c2}
B.R. Barmish, "A Generalized of Kharitonov's Four Polynomial
concept for Robustness Stability with Linearly dependent
coeficient Perturbations" {\it in IEEE Transactions on Automatic
Control}, Vol. 34, pp. 157--165, 1989.

\bibitem{c3}
Z. Gong, C. Wen and D.P. Mital, "Decentralized Robust Controller
Design for a Class of Interconnected Uncertain Systems: With
Unknown Bound of Uncertainty", {\it IEEE Transactions on Automatic
Control}, Vol. 41, No. 6, pp. 850-854, 1996.

\bibitem{c9}
F. Liao, J.L. Wang,G-H. Yang, "Reliable Robust Flight Tracking
Control: An LMI Approach".{\it IEEE Trans. on control Systems
Techn.}, Vol. 10, No. 1, pp.76-89, 2002.

\bibitem{c4}
H. Mukaidini, Y. Takato, Y. Tanaka and K. Mizukami, "The
Guaranteed Cost Control for Uncertain Large-scale Interconnected
Systems", {\it in 15th Triennial World Congress}, Barcelona,
Spain, 2002.

\bibitem{c8}
Pujol G., "Contribution on reliable control for uncertain
interconnected systems" {\it PhD Thesis Department of Applied
Mathematics III}, Technical University of Catalonia, 2004.

\bibitem{c5}
H.D. Tuan, P. Apkarian and T.Q. Nguyen, "Robust Filtering for
Uncertain Nonlinear Parameterized Plants" {\it IEEE Transactions
on Signal Processing}, Vol. 51, No. 7, pp. 1816--1824, 2003.

\bibitem{c6}
S. Xie, L. Xie, Y. Wang and H. Zhang, "Decentralized Guaranteed
cost Control of a Class of Large-Scale Interconnected Systems",
{\it Proc. $38^{th}$ Conference on Decision and Control, Phoenix,
Arizona, USA}, pp. 3297--3302, 1999.

\bibitem{c7}
G.-H. Yang, J.L. Wang and Y.C. Soh, "Reliable Guaranteed Cost
Control for Uncertain Nonlinear Systems", {\it IEEE Transactions
on Automatic Control}, Vol.45, No.11, pp. 2188--2192, 2000.

\bibitem{Sil80}
D.D. Siljak,"Reliable control using redundant controllers", {\it
IEEE Transactionss on Automatic Control}, Vol.31,
pp.303--329,1980.

\bibitem{Yang98}
G.-H. Yang, S.-Y. Zhang, J. Lam and J. Wang, "Reliable control
using redundant controllers", {\it IEEE Transactions on Automatic
Control}, Vo.43, pp.1588--1593, 1998.


\end{thebibliography}

\end{document}